# Über die Äquisingularität von normalen Flächensingularitäten

Achim Hennings ♦

**Abstract**
We give a description of the equisingularity of a family of normal surface singularities by some numerical invariants belonging to them. We consider two notions of equisingularity, Whitney regularity and a more restrictive condition using the Nash transformation.

## 0 Einleitung

Sei $X \to D$ eine flache Familie von normalen Flächensingularitäten über der Kreisscheibe $D \subseteq \mathbb{C}$ um 0, deren kritische Menge glatt über $D$ liegt. Sei $X \subseteq B \times D$, $B \subseteq \mathbb{C}^N$, ein geeigneter Repräsentant, so dass $\{0\} \times D$ die kritische Menge der Projektion ist. Zu $X$ bildet man die relative Nash-Transformation, ebenso für den Schnitt $X_H$ mit einer allgemeinen Hyperebene durch $\{0\} \times D \subseteq \mathbb{C}^N \times D$.

Im ersten Teil dieser Arbeit betrachten wir die folgende Äquisingularitätsbedingung:[1]

(1) Die relative Nash-Transformation von $X$ hat konstante Faserdimension.
(2) Die relative Nash-Transformation von $X_H$ hat konstante Faserdimension.
(3) Die Multiplizität der Fasern von $X$ ist konstant.

Die Bedingungen (2), (3) beschreiben gerade die Whitney-Regularität von $X_H$, da man in Dimension eins die Nash-Transformation auch durch den Konormalenraum ersetzen darf. Bei (1) ist dies jedoch nicht mehr der Fall, so dass die Bedingung insgesamt stärker als die Whitney-Regularität ist.

Die so definierte Äquisingularität soll durch numerische Invarianten der Fasern charakterisiert werden. Es ist bekannt, dass die Whitney-Regularität von Kurven durch die Milnor-Zahl und die Multiplizität ausgedrückt werden kann. Zusätzlich betrachten wir die erste Chern-Klasse des Nash-Bündels, die als lokale Klasse repräsentierbar ist. Quadriert, ergibt sich als Invariante eine rationale Zahl. Nach einem Satz von M. Kato kann diese Zahl als verallgemeinerte Multiplizität[2] des definierenden Ideals der Nash-Transformation angesehen werden und ähnlich wie im Hyperflächen-Fall ist das Spezialisierungsprinzip für ganze Abhängigkeit von B. Teissier anwendbar.

Im zweiten Teil der Arbeit untersuchen wir die Whitney-Regularität von $X$. Nach einem Ergebnis von J. Kollar und N. Shepherd-Barron besitzt $X$ eine sehr schwache simultane Auflösung, wenn die Schnittzahl $K_t^2$ des kanonischen Bündels der minimalen Auflösung der Fasern konstant ist. Nimmt man die Eulercharakteristik der minimalen Auflösung, sowie

---

♦ Universität Siegen, Hölderlinstraße 3, 57068 Siegen

[1] Äquisingularitätsbedingungen dieser Art wurden unter allgemeinen Gesichtspunkten von J.-P. Henry und M. Merle [HM] untersucht.
[2] In der neueren Literatur wird auch die höherdimensionale Verallgemeinerung betrachtet und teilweise als $\varepsilon$−Multiplizität bezeichnet.

Milnor-Zahl und Multiplizität des allgemeinen Hyperebenenschnitts hinzu, erhält man eine Charakterisierung der Whitney-Regularität.

**1 Erste Chern-Klasse des Nash-Bündels**

Sei $(X, 0) \subseteq (\mathbb{C}^N, 0)$ eine normale Singularität der Dimension 2, in einer guten (Kugel-) Umgebung $B$ definiert durch das Ideal $I$. Die Nash-Transformation ist geometrisch der Abschluss des Bildes der Abbildung $X - \{0\} \to X \times G, x \mapsto (x, T_x)$, in das Grassmann-Bündel, die glatten Punkten ihren Tangentialraum $\subseteq \mathbb{C}^N$ zuordnet. Man erhält eine eigentliche Modifikation $\nu: \hat{X} \to X$ und ein Vektorraumbündel $T \to \hat{X}$, welches das Tangentialbündel fortsetzt. Die Nash-Transformierte kann als Aufblasung eines Ideals dargestellt werden. Sei $I_1 = (g_1, \dots, g_{N-2})$ das Ideal eines vollständigen Durchschnitts $X_1 \supseteq X$, der außerhalb einer Hyperfläche $Z$, die $X$ nicht enthält, mit $X$ übereinstimmt (vgl. Anhang 2). Seien $\delta_1, \dots, \delta_m$ die $m = \binom{N}{N-2}$ maximalen Minoren der Jacobi-Matrix $(g_{x_1}, \dots, g_{x_N})$ und $J$ das Ideal $(\delta_1, \dots, \delta_m)\mathcal{O}_X$. Dieses Ideal ist außerhalb des Nullpunkts ein Hauptideal. Für die glatten Punkte $x \in X - Z$ beschreiben die Zuordnungen

$$x \mapsto [g_{x_1}(x), \dots, g_{x_N}(x)] \mapsto [\delta_1(x), \dots, \delta_m(x)]$$

die Abbildungen

$$x \mapsto (\mathbb{C}^N/T_x)^* \mapsto \Lambda^{N-2}(\mathbb{C}^N/T_x)^*,$$

die eine Einbettung von $\hat{X}$ in $X \times \mathbb{P}(\Lambda^{N-2}(\mathbb{C}^N)^*)$ liefern. Also ist $\hat{X}$ die Aufblasung von $J \subseteq \mathcal{O}_X$, denn sie entsteht ebenso wie die Nash-Transformierte durch Abschließung. Das inversible Ideal $J\mathcal{O}_{\hat{X}}$ ist die Garbe der Schnitte des Geradenbündels

$$L := \Lambda^{N-2}(\mathbb{C}^N/T)^{**} = \Lambda^{N-2}(\mathbb{C}^N/T) = \Lambda^2 T^*.$$

Unter der exzeptionellen Menge verstehen wir hier das kompakte Urbild $\hat{E} := \nu^{-1}(0)$. Der Divisor von $J\mathcal{O}_{\hat{X}}$ (der üblicherweise exzeptioneller Divisor genannt wird) kann dagegen noch zusätzlich nicht kompakte Komponenten haben (die außerhalb $\hat{E}$ isomorph abgebildet werden).

Sei $\rho: \tilde{X} \to \hat{X}$ eine Singularitätenauflösung und $E := \rho^{-1}(\hat{E})$ die exzeptionelle Menge von $\nu\rho$. Dabei seien $\tilde{X}, \hat{X}$ Repräsentanten ohne Rand, in denen die exzeptionelle Menge ein starker Deformationsretrakt ist[3]. Die entsprechenden berandeten Räume bezeichnen wir mit $\overline{\tilde{X}}, \overline{\hat{X}}$ und den Rand mit $\partial \overline{X}$. Wegen der Nicht-Ausartung der Schnittform auf $\tilde{X}$ ist in der exakten Sequenz (mit komplexen Koeffizienten)

$$H^2_E(\tilde{X}) \xrightarrow{\alpha} H^2(\tilde{X}) \xrightarrow{\beta} H^2(\tilde{X} - E)$$

die Abbildung $\alpha$ surjektiv (sogar bijektiv) und $\beta$ null. Also ist auch die Komposition $H^2(\hat{X}) \to H^2(\tilde{X}) \to H^2(\tilde{X} - E)$ null und $H^2_{\hat{E}}(\hat{X}) \to H^2(\hat{X})$ surjektiv. Daher hat die erste Chern-Klasse $c_1(T^*) \in H^2(\hat{X})$ des dualen Nash-Bündels ein Urbild $\tilde{c} \in H^2_{\hat{E}}(\hat{X}) \cong H^2(\overline{\hat{X}}, \partial \overline{X}) \cong$ $H^2_c(\hat{X})$.[4] Dann ist das Produkt $\tilde{c}^2 \in H^4_{\hat{E}}(\hat{X})) \cong H^4_c(\hat{X})$. wohldefiniert und wir betrachten dessen Grad

$$\deg(\tilde{c}^2) \coloneqq \deg(\tilde{c}^2 \cap [\hat{X}]) = \int_{\hat{X}} \tilde{c}^2.$$

---

[3] Die Aussage für $\hat{X}$ folgt unmittelbar aus der für $\tilde{X}$.

Man kann auch $\rho^*\tilde{c} \in H^2_c(\tilde{X})$ bilden. Dies ist das Urbild für $c_1(\rho^*T^*) \in H^2(\tilde{X})$, und es gilt $\int_{\hat{X}} \tilde{c}^2 = \int_{\tilde{X}} (\rho^*\tilde{c})^2$.

Bezeichnung: Das Bündel $L$ bezeichnen wir auch mit $\hat{K}$, da es eine Fortsetzung des kanonischen Bündels ist, und schreiben $\hat{K}^2$ für die soeben definierte Zahl.

Bemerkung: Seien $\hat{X}_1 \subseteq \hat{X}$ zwei Repräsentanten, so dass die Restriktion $H^2(\hat{X}) \to H^2(\hat{X}_1)$ isomorph ist, und $\tilde{c} \in H^2_c(\hat{X})$ bzw. $\tilde{c}_1 \in H^2_c(\hat{X}_1)$ Urbilder von $c \in H^2(\hat{X})$ bzw. $c|\hat{X}_1 \in H^2(\hat{X}_1)$. Dann ist $\deg \tilde{c}^2 = \deg \tilde{c}_1^2$. Denn das Diagramm

$$\begin{array}{ccc} H^2(\hat{X}) & \to & H^2(\hat{X}_1) \\ \uparrow & & \uparrow \\ H^2_c(\hat{X}) & \leftarrow & H^2_c(\hat{X}_1) \end{array}$$

der kanonischen Abbildungen ist kommutativ. Also wird auch $\tilde{c}_1$ auf $c$ abgebildet. Das Produkt $\tilde{c}^2$ hängt aber nur von dem Bild $c$ ab. Die hier gemachte Voraussetzung ist erfüllt, wenn $\hat{X}_1, \hat{X}$ Repräsentanten in schrumpfenden guten Umgebungen sind, kann aber bei Deformation von $X$ verloren gehen.

Die bisherigen Überlegungen gelten auch für allgemeine Aufblasungen an Stelle der Nash-Transformation, wenn das Zentrum auf $X^* \coloneqq X - \{0\}$ eine Hyperfläche ist. Sei also $J \subseteq \mathcal{O}_X$ ein Ideal, das auf $X^*$ inversibel ist. In diesem Fall ist $\hat{X}$ die Aufblasung des Ideals, $L$ das Geradenbündel zu $J\mathcal{O}_{\hat{X}}$ und $\tilde{c}$ ein Urbild von $c_1(L) \in H^2(\hat{X})$ und wir betrachten $\deg(\tilde{c}^2)$. Die topologische Deutung dieser rationalen Zahl kann man mit der Schnittform konkretisieren: Sei $E = E_1 \cup \ldots \cup E_r$ die Zerlegung der exzeptionellen Menge der Auflösung $\rho: \tilde{X} \to \hat{X}$ in irreduzible Komponenten. Sei $l$ das Urbild von $L$ auf $\tilde{X}$. Dann gibt es einen eindeutigen rationalen Zyklus $z = \sum_{i=1}^r a_i E_i$ mit $l \cdot E_i = z \cdot E_i$ für alle $i$. Dieser stellt das Bild von $\rho^*\tilde{c}$ in der Homologie dar. Es gilt daher

$$\deg(\tilde{c}^2) = z \cdot z.$$

Nach einem Satz von M. Kato [K] (vgl. [W]) kann $l \cdot l \coloneqq z \cdot z$ auch als verallgemeinerte Multiplizität gedeutet werden. Dort wird nämlich gezeigt, dass,

$$\dim_{\mathbb{C}} \frac{\Gamma(X^*, \mathcal{O}(l^n))}{\Gamma(\tilde{X}, \mathcal{O}(l^n))} + h^1(\mathcal{O}(l^n)) = (-l^2)\frac{n^2}{2} + (l \cdot K)\frac{n}{2} + b(n)$$

wobei $K$ das kanonische Bündel, $l \cdot K \coloneqq z \cdot K$ und $b$ eine beschränkte Funktion ist. Da $\dim R^1\rho_*\mathcal{O}(l^n)$ konstant ist (wie die Faktorisierung von $\rho$ über die Normalisierung von $\hat{X}$ zeigt) und $h^1(\mathcal{O}(L^n) \otimes \rho_* \mathcal{O}_{\tilde{X}}) = 0$, $n \gg 0$, ist $h^1(\mathcal{O}(l^n))$ beschränkt und trägt nichts zu $l^2$ bei. Aus diesem Grund schreiben wir $-l^2 = -\deg(\tilde{c}^2)$ auch in der Form[5]

$$e(J, 0) = e(J).$$

---

[4] $\tilde{c}$ kann durch eine glatte geschlossene 2-Form repräsentiert werden. Denn $c_1(T^*)$ wird durch eine 2-Form $\varphi$ repräsentiert und die Beschränkung auf den *glatten* Raum $\tilde{X} - E$ ist die Ableitung einer 1-Form $\psi$. Mit einer Abschneidefunktion $\sigma$ bildet man dann $\varphi - d((1-\sigma)\psi)$.

Wenn $J' \subseteq J$ ein zweites Ideal ist, das außerhalb des Nullpunkts mit $J$ übereinstimmt, dann ist die Differenz $e(J') - e(J)$ das zweifache des Leitkoeffizienten zu der polynomialen Funktion[6]

$$\dim_{\mathbb{C}}(J^n/J'^n), n \gg 0.$$

Dies ist die in [R] und [KT1] untersuchte Multiplizität von Paaren und es wird dort [R, th. 2.1], [KT1, prop. 5.3] gezeigt, dass sie genau dann Null ist, wenn $J' \subseteq J$ eine Reduktion ist.

## 2 Faserdimension der relativen Nash-Transformation

Sei jetzt $\tau: X \to D$ eine Deformation von zweidimensionalen normalen Singularitäten mit einer Kreisscheibe $D$ um 0 als Basis und kritischer Menge $C$. Dabei sei $X \subseteq B \times D$ und $B \subseteq \mathbb{C}^N$ eine gute Umgebung für die spezielle Faser $X_0$. Sei $v: \hat{X} \to D$ die relative Nash-Transformation, $\hat{E} = v^{-1}(C)$ die exzeptionelle Menge und $T$ das fortgesetzte relative Tangentialbündel. Die relative Nash-Transformation ist die Aufblasung eines Ideals, das außerhalb von $C$ inversibel ist.

Wir nehmen nun an, dass $\hat{X}_{t,red}$ für kleine $t$ die Nash-Transformation der Faser, d.h. die strikte Transformierte von $X_t$ in $\hat{X}$ ist. Dies ist genau dann der Fall, wenn $\hat{E}_0 = \hat{E} \cap \hat{X}_0$ keine Komponente von $\hat{X}_0$ enthält, also wenn $\dim \hat{E}_0 \leq 1$. Aufgrund der Wahl des Repräsentanten ist die Nullfaser $\hat{X}_0$ auf die exzeptionelle Menge $\hat{E}_0$ kontrahierbar. Daher ergibt

$$H^2(\hat{X}) \to H^2(\hat{X}_0) \to H^2(\hat{X}_0 - \hat{E}_0) \cong H^2(\partial \bar{X}_0)$$

wieder die Nullabbildung. Da $\partial \bar{X} \to D$ trivial ist, gilt $H^2(\partial \bar{X}_t) \cong H^2(\partial \bar{X})$ für alle $t$. Also ist auch

$$H^2(\overline{\hat{X}}) \to H^2(\partial \bar{X})$$

null und daher

$$H^2(\overline{\hat{X}}, \partial \bar{X}) \to H^2(\overline{\hat{X}}) \cong H^2(\hat{X})$$

surjektiv. Sei $\tilde{c} \in H^2(\overline{\hat{X}}, \partial \bar{X})$ ein Urbild von $c_1(T^*)$.[7]

Wir nehmen jetzt zusätzlich an, dass $C$ glatt über $D$ und die Deformation Whitney-regulär ist. Dann bleibt $B$ eine gute Umgebung für die Fasern $X_t$, $t$ nahe 0, und es gilt (nicht nur für $t = 0$)

---

[5] In neueren Arbeiten wird eine entsprechende Begriffsbildung (z.T. unter der Bezeichnung $\varepsilon$ −Multiplizität) auch in höherer Dimension untersucht, hat dann allerdings keine so günstigen Eigenschaften mehr.

[6] Zunächst $\dim_{\mathbb{C}}(\overline{J^n}/\overline{J'^n})$, aber der ganze Abschluss ändert den Leitkoeffizienten nicht.

[7] $\tilde{c}$ kann wieder durch eine glatte geschlossene 2-Form mit eigentlichem Träger repräsentiert werden.

$$H^2(X_t - C_t) \cong H^2(\partial \bar{X}_t),$$

$$H^2_{\hat{E}_t}(\hat{X}_t) \cong H^2(\overline{\hat{X}}_t, \partial \bar{X}_t).$$

Also ist $\tilde{c}(t) := \tilde{c}|\hat{X}_t$ eine geeignete Wahl für die Faser $X_t$ im Sinne von oben. Der Grad

$$\deg(\tilde{c}(t)^2) = \int_{\hat{X}_t} \tilde{c}(t)^2$$

ist konstant. Denn der reduzierte Zyklus $[\hat{X}_t] - [\hat{X}_0]$ ist nullhomolog.[8]

Diese Betrachtungen können wieder auf allgemeinere Aufblasungen übertragen werden. Sei $J_1 \subseteq \mathcal{O}_X$ ein kohärentes Ideal und $C_1$ endlich über $D$, so dass $J_1$ außerhalb von $C_1$ ein Hauptideal ist und $J_1(t) := J_1 \mathcal{O}_{X_t}$ nicht Null ist. Dann ist $e(J_1(0), 0)$ definiert. Ebenso definiert man

$$e(J_1(t)) := \sum_{x \in C_{1,t}} e(J_1(t), x).$$

Da eine Auflösung $\tilde{X}$ von $X$, die über die Aufblasung $\hat{X} \to X$ von $J_1$ faktorisiert, für $t \neq 0$ nahe 0 eine schwache simultane Auflösung ist, ist die oben definierte Selbstschnittzahl von $J_1 \mathcal{O}_{\tilde{X}_t}$ und damit $e(J_1(t))$ konstant für kleine $t \neq 0$.

Wir können nun die Aussage dieser Arbeit formulieren. Sei wie vorher $\tau: X \to D$ eine Deformation von normalen Flächensingularitäten über der Kreisscheibe $D \subseteq \mathbb{C}$ um 0 mit dem Repräsentanten $X \subseteq B \times D$, $B \subseteq \mathbb{C}^N$ gute Umgebung für die zentrale Faser, und wir nehmen an, dass die kritische Menge $C = \{0\} \times D$ ist. Sei $X_H$ der Durchschnitt mit einer allgemeinen Hyperebene $H = H' \times D$, $H' \subseteq \mathbb{C}^N$, die $C$ enthält. $X$, $X_H$ und $X_{H,t}$ sind als Cohen-Macaulay-Räume reduziert, da ihre Singularitätenmenge niederdimensional ist. Wir schreiben $\mu^{(2)}(X_t, 0)$ für die Milnor-Zahl von $X_{H,t}$, wenn $H'$ (für dieses $t$) allgemein ist. Da $\mu^{(2)}(X_t, 0)$ und $\mu(X_{H,t}, 0)$ nur isolierte Sprungstellen haben, gilt für allgemeine $H$: $\mu^{(2)}(X_t, 0) = \mu(X_{H,t}, 0)$ für alle $t$ nahe 0. Eine ähnliche Bemerkung gilt für die Multiplizität: $m(X_t, 0) = m(X_{H,t}, 0)$ für $t$ nahe 0, falls $H$ allgemein ist.

**Satz 1**: Es sind äquivalent:

a) (1) Die relative Nash-Transformation $\hat{X}$ von $X$ hat Faserdimension $\leq 1$.
   (2) Die relative Nash-Transformation $\hat{X}_H$ von $X_H$ hat Faserdimension $0$.
   (3) Die Multiplizität der Fasern von $X_H$ (oder $X$) ist konstant längs $C$.
b) Die Multiplizität $m(X_t, 0)$ und die Milnor-Zahl $\mu^{(2)}(X_t, 0)$ sind konstant längs $C$ und der Grad der quadrierten ersten Chern-Klasse $\hat{K}_t^2$ des Nash-Bündels von $(X_t, 0)$ ist konstant.

Beweis: a) $\Rightarrow$ b). Aus den Voraussetzungen folgt, dass die Deformation Whitney-regulär ist (vgl. Anhang 1). Dies gilt dann auch für $X_H$, insbesondere ist $C$ die kritische Menge, und die

---

[8] Sei $\sigma$ eine glatte 2-Form mit kompaktem Träger um $0 \in \mathbb{C}$ und $\int_{\mathbb{C}} \sigma = 1$. Dann gilt für glatte geschlossene 4-Formen $\alpha$ auf $\hat{X}$ mit eigentlichem Träger $\int_{\hat{X}_t} \alpha = \int_{\hat{X}} \alpha \wedge ((\tau - t)v)^* \sigma$. Dieses Integral ist unabhängig von $\sigma$ und $t$.

Milnor-Zahl und die Multiplizität sind konstant [BGG, th. III.3], vgl. [BG, p. 264]. Oben wurde erläutert, dass auch $\widehat{K}_t^2$ konstant ist.

b) ⇒ a). Wir zeigen (1): Sei $J \subseteq \mathcal{O}_X$ ein definierendes Ideal der Nash-Transformation (inversibel außerhalb $C$ und $J(0) \neq 0$). Nach Annahme ist $C = \{0\} \times D$ glatt über $D$ und $e(J(t))$ ist konstant für $t$ nahe $0$. Sei $J_1 \subseteq J$ erzeugt von 2 Elementen, so dass $J_1(0) \subseteq J(0)$ eine Reduktion ist. Sei $Z$ der Träger von $J/J_1$. Da ein Hauptideal keine echte Reduktion haben kann, stimmen $J_1(0) \subseteq J(0)$ auf $X_0 - \{0\}$ überein. Dann sind auch $J_1 \subseteq J$ in einer Umgebung von $X_0 - \{0\}$ in $X$ gleich.[9] Daher kann $Z$ keine 1–kodimensionale Komponente haben weil diese wegen der Halbstetigkeit der Faserdimension auch $X_0 - \{0\}$ treffen würde. Die Ideale $J_1 \subseteq J$ unterscheiden sich also nur in Kodimension 2, d.h. nur in endlichen Teilmengen der Fasern $X_t$. In diesem Fall ist klar, dass

$$e(J_1(t)) \geq e(J(t)).$$

Nach dem Lemma unten ist die Halbstetigkeit $e(J_1(0)) \geq e(J_1(t))$ gültig. Daraus folgt mit der Ganz-Abhängigkeit

$$e(J(0)) = e(J_1(0)) \geq e(J_1(t)) \geq e(J(t)) = e(J(0)).$$

Also ist $e(J_1(t)) = e(J(t))$ für alle $t$. Insbesondere folgt $e(J_1(t), x) = 0$ für $x \in X_t - C$ und die Inklusion $J_1 \subseteq J$ ist außerhalb von $C$ eine Gleichheit. Weiter gilt $e(J_1(t), 0) = e(J(t), 0)$ für alle $t$. Wie in Abschnitt 1 schließt man daraus zunächst $\overline{J_1(t)}_0 = \overline{J(t)}_0$ für $t \neq 0$. Der exzeptionelle Divisor der Aufblasung von $J$ ist flach über $D$ für kleine $t \neq 0$.[10] Denn die über $D$ nicht flachen Punkte liegen im Urbild von $C$ und ihr Bild ist eine echte analytische Teilmenge in $D$ (Satz von J. Frisch, vgl. [BS]; im algebraischen Fall Satz über generische Flachheit). Daher ist die Koordinate $\tau - t$ superfiziell (vom Grad 0) für das Ideal $J$ und man schließt $\overline{J}_{1,(0,t)} = \overline{J}_{(0,t)}$ für $t \neq 0$. Da $J_1$ von 2 Elementen erzeugt wird, hat $\overline{J}_1$ nur assoziierte Primideale der Höhe $\leq 2$. (Ein Element ist ganz über $J_1$, wenn die Ordnung längs jeder Komponente des exzeptionellen Divisors der normalisierten Aufblasung von $J_1$ mindestens so groß wie die von $J_1$ ist. Weil diese Komponenten Dimension 2 haben, liegt keine in der speziellen Faser.- Alternativ kann man das algebraische Analogon des Riemannschen Hebbarkeitssatzes auf das affine Spektrum des normalisierten Rees-Rings von $J_{1,(0,0)}$ und das Komplement der speziellen Faser anwenden.) Daraus folgt $\overline{J}_1 = \overline{J}$. Damit hat die Nash-Transformation konstante Faserdimension.

Nun hat $X_H$ die kritische Menge $C$, wenn $H'$ kein Element der speziellen Faser von $\widehat{X}$ enthält (vgl. Anhang 1). Aus der konstanten Milnor-Zahl und Multiplizität folgt dann nach [BGG, th. III.3] die Whitney-Regularität von $X_H$, also (2), (3).

**Lemma:** (Halbstetigkeit der verallgemeinerten Multiplizität)

---

[9] Bei $x \in X_0 - \{0\}$ sei $J_x = (a)$ und $J_1(0)_x = (b|_{t=0})$ mit $b = ca \in J_{1x}$, $b|_{t=0} = a|_{t=0}$. Dann ist $c(x) \neq 0$ und $c$ eine Einheit.

[10] Es genügt, dass die Faserdimension konstant ist.

Sei $J_1 \subseteq \mathcal{O}_X$ ein kohärentes Ideal, inversibel außerhalb der über $D$ endlichen Menge $C_1$ mit $J_1(0)$ nicht null. Dann ist $e(J_1(t)) \leq e(J_1(0))$.

Beweis: Seien $j: X - C_1 \to X$, $j_t: X_t - C_{1,t} \to X_t$ die Inklusionen und $\widetilde{J_1^n} := j_* j^{-1} J_1^n$, $\widetilde{J_1^n(t)} := j_{t*} j_t^{-1} J_1^n(t)$.

Behauptung: $\dim \widetilde{J_1^n(t)} / J_1^n(t) \leq \dim \widetilde{J_1^n(0)} / J_1^n(0)$.

Beweis: Die Faser über $t \in D$ wird definiert durch die Funktion $\tau - t$. Sei in diesem Beweis $J := \widetilde{J_1^n}, I := J_1^n$.

(1) Nach Definition ist $\tau - t$ kein Nullteiler in $\mathcal{O}_X / J$. Daher gilt $J \otimes \mathcal{O}_{X_t} = J \mathcal{O}_{X_t}$.

(2) $M := \tau_*(J/I)$ ist ein kohärenter $\mathcal{O}_D$ −Modul und daher torsionsfrei für $t \neq 0$. Also ist für $t \neq 0$ die Abbildung $I \otimes \mathcal{O}_{X_t} \to J \otimes \mathcal{O}_{X_t}$ injektiv. Die Dimension von

$$M(t) = \tau_*(J \otimes \mathcal{O}_{X_t} / Bild(I \otimes \mathcal{O}_{X_t})),$$

ist für $t \neq 0$ konstant, und $\dim M(0) \geq \dim M(t)$.

(3) Wir zeigen $\widetilde{J_1^n}(t) = \widetilde{J_1^n(t)}$ für $t \neq 0$. Wir betrachten $N := J_1^n$ nur als Modul und bezeichnen mit Klammern $(t)$ das Tensorprodukt $\otimes \mathcal{O}_{X_t}$. Zu zeigen ist $N^{**}(t) = N(t)^{**}$. Aus der exakten Sequenz $0 \to \mathcal{O}_X \to \mathcal{O}_X \to \mathcal{O}_{X_t} \to 0$ folgt

$$0 \to N^* \to N^* \to N(t)^* \to \underline{Ext}^1_{\mathcal{O}_X}(N, \mathcal{O}_X) \to \underline{Ext}^1_{\mathcal{O}_X}(N, \mathcal{O}_X).$$

Der Modul $\underline{Ext}^1_{\mathcal{O}_X}(N, \mathcal{O}_X)$ hat Träger in $C_1$ und ist endlich über $D$, also für $t \neq 0$ torsionsfrei. Daraus folgt $N^*(t) \cong N(t)^*$ und analog $N^{**}(t) \cong N^*(t)^* \cong N(t)^{**}$.

Für $t = 0$ gilt $\widetilde{J_1^n}(t) \subseteq \widetilde{J_1^n(t)}$.

Dann ist wegen (1)

$$M(0) \cong \widetilde{J_1^n}(0) / J_1^n(0) \subseteq \widetilde{J_1^n(0)} / J_1^n(0).$$

Daraus folgt die Behauptung.

Mit dem oben benutzten Satz von M. Kato folgt hieraus die Halbstetigkeit von $e(J_1(t))$, wenn man beachtet, dass die Koeffizienten des Polynoms und die Fehlerschranken aufgrund simultaner Auflösung über $t \neq 0$ dort konstant sind.

**3 Invarianten für die Whitney-Regularität**

Wie in Abschnitt 2 sei $X \subseteq B \times D$ eine flache Familie von normalen Flächensingularitäten mit kritischer Menge $C = \{0\} \times D$ und sei $B \subseteq \mathbb{C}^N$ gute Umgebung für die spezielle Faser $X_0$. Wir wählen auch gute Umgebungen $B_t \subseteq B$ für die Nachbarfasern $X_t$, $t$ nahe 0. Um die Ergebnisse aus [KSB] anzuwenden, setzen wir voraus, dass $X$ zu einer Deformation von projektiven algebraischen Flächen erweitert werden kann. Nach [A], [E], [L1] stellt dies keine

Einschränkung dar. Sei $\pi_t: \tilde{X}_t \to X_t$ die minimale Auflösung und $E_t = \pi_t^{-1}(0)$ die exzeptionelle Menge. Das kanonische Bündel $K_t$ von $\tilde{X}_t$ kann in der Nähe von $E_t$ (d.h. über $X_t \cap B_t$) durch einen numerisch äquivalenten exzeptionellen Zyklus ersetzt werden. Auf diese Weise ist die Schnittzahl $K_t^2$ definiert und hängt nur von der Singularität $(X_t, 0)$ ab (vgl. Abschnitt 1). Nach [KSB, Beweis von Th. 2.22] gilt: Wenn die Schnittzahl $K_t^2$ konstant in $t$ ist, dann gibt es (nach endlichem Basiswechsel) eine simultane RDP-Auflösung[11] und damit eine sehr schwache simultane Auflösung $\pi: \tilde{X} \to X$, wo die Fasern $\tilde{X}_t$ die minimale Auflösung von $X_t$ sind. Dann ist $\tilde{X} \to D$ topologisch trivial, insbesondere $\chi(\tilde{X}_t) = \chi(\tilde{X}_0) = \chi(E_0)$. Es gilt, da Umgebungsränder ungerade Dimension haben und keinen Beitrag zur Eulercharakteristik liefern:

$$\chi(X_t) = \chi(X_t \cap B_t) + \chi(X_t - B_t),$$

$$\chi(\tilde{X}_t) = \chi(\tilde{X}_t \cap \pi^{-1} B_t) + \chi(X_t - B_t) = \chi(E_t) + \chi(X_t - B_t).$$

Also folgt $\chi(X_t) = \chi(\tilde{X}_t) - \chi(E_t) + 1$. Wenn man $\chi(E_t)$ als konstant annimmt, sind also auch $\chi(X_t) = 1$ und $\chi(X_t - B_t) = 0$ konstant.

Sei $N(X/D)$ die relative Nash-Transformation mit dem fortgesetzten Tangentialbündel $T$ und $C(X/D)$ der relative Konormalenraum. Die entsprechenden absoluten Räume der Fasern $N(X_t)$ und $C(X_t)$ können als strikte Transformierte in die relativen Räume eingebettet werden. In diesem Sinne gilt

$$[N(X/D)_t] = [N(X_t)], t \neq 0,$$

$$[N(X/D)_0] = [N(X_0)] + \hat{E}_v, [C(X/D)_0] = [C(X_0)] + F_v,$$

wobei $\hat{E}_v$ bzw. $F_v$ der Anteil der Zyklen über $0$ ist.

Sei $s$ ein differenzierbarer Schnitt von $T$ über $N(X/D)$, der in der Nähe des Randes $\partial \bar{X} = \bar{X} \cap \partial B \times D$ radial ist. Man erhält eine lokalisierte Chern-Klasse $c_2(T, s)$, indem man das Urbild der Thom-Klasse $\tau \in H^4(T, T - 0)$ unter $s$ betrachtet. Alternativ kann man einen differenzierbaren Zusammenhang mit $\nabla s = 0$ am Rand zugrunde legen. Es gibt eine Darstellung der Klasse durch eine Differentialform mit eigentlichem Träger. Da man (für jedes feste $t$) annehmen kann, dass $s$ auf $X_t - B_t$ nur isolierte Nullstellen hat und radial bei $X_t \cap \partial B_t$ ist, gilt nach Definition der Euler-Obstruktion und nach der Indexformel für Vektorfelder

$$\deg(c_2(T, s) \cap N(X_t)) = \chi(X_t - B_t) + Eu(X_t, 0).$$

Nach [LT, prop. 6.2.1] ist die Euler-Obstruktion $Eu(X_t, 0)$ die Eulercharakteristik der Milnor-Faser einer allgemeinen Linearform auf $X_t$, also $1 - \mu^{(2)}(X_t, 0)$, da nach [BGG] und [BG] die Milnor-Faser zusammenhängend ist.

---

[11] Minimale Auflösung mit rationalen Doppelpunkten

Wir nehmen an, dass $\mu^{(2)}(X_t, 0)$ konstant ist. Dann ist also $\deg(c_2(T, s) \cap N(X_t))$ konstant. Da $\deg(c_2(T, s) \cap N(X/D)_t)$ generell konstant ist[12], folgt

$$\deg(c_2(T, s) \cap \hat{E}_v) = \deg(c_2(T) \cap \hat{E}_v) = 0.$$

Sei $L$ das duale universelle Bündel auf dem projektiven Bündel $\mathbb{P}((\mathbb{C}^N/T)^*) \to N(X/D)$. Dies ist das Urbild des Hyperebenenbündels auf $C(X/D)$, das wir auch mit $L$ bezeichnen, unter der surjektiven Abbildung

$$\mathbb{P}((\mathbb{C}^N/T)^*) \to C(X/D),$$

die man durch Vernachlässigung der Punkte in $N(X/D)$ erhält.

Nach der Projektionsformel wird $\mathbb{P}((\mathbb{C}^N/T)^*)|\hat{E}_v$ auf $F_v$ abgebildet. Daher ist

$$(-1)^2 \deg(c_2(T) \cap \hat{E}_v) = \deg(L^{N-1} \cap [\mathbb{P}((\mathbb{C}^N/T)^*)|\hat{E}_v]) = \deg(L^{N-1} \cap F_v).$$

Da dies null ist, muss der Schnitt von $F_v$ mit $N-1$ Hyperebenen in $\mathbb{P}^{N-1}$ leer sein, d.h. die spezielle Faser von $C(X/D)$ hat Dimension $\leq N - 2$. Dies charakterisiert zusammen mit der Whitney-Regularität des allgemeinen Hyperebenenschnitts $X_H$, die durch Milnor-Zahl und Multiplizität bestimmt wird, die Whitney-Regularität von $X$ (Anhang 1 und [BGG]).

Da hieraus die schwache simultane Auflösbarkeit folgt [KSB, Th. 2.10] [13], erhält man:

**Satz 2:** Notwendig und hinreichend für die Whitney-Regularität von $X$ ist, dass die Invarianten

$$K_t^2, \chi(E_t), \mu^{(2)}(X_t, 0), m(X_t, 0)$$

konstant sind

Bemerkung: Mit Lemma 2.20 in [KSB] und ähnlichen (aber elementaren) Argumenten wie in Th. 2.22 in [KSB] erhält man auch einen alternativen Beweis für den Satz 1 aus Abschnitt 2.

**Anhang 1: Beschreibungen der Whitney-Regularität**

Sei $X \subseteq U \times D$, $U \subseteq \mathbb{C}^N$ offen, eine Familie von 2-dimensionalen isolierten Singularitäten über $0 \in D \subseteq \mathbb{C}$ mit kritischer Menge $C \supseteq Y = \{0\} \times D$ für die Projektion $p: X \to D$. Da wir nur reduzierte Räume betrachten wollen, setzen wir der Einfachheit halber voraus, dass die Familie flach ist und die Singularitäten normal sind. Sei $\kappa: C(X/D) \to X$ der relative Konormalenraum und $B_Y C(X/D)$ die Aufblasung längs $\tilde{Y} = \kappa^{-1}(Y)$. Dies ist der Abschluss von

$$\{(x, h, l) | x = (x_1, \ldots, x_N, t) \in X - C, T_x \subseteq h, l = [x_1, \ldots, x_N]\} \subseteq X \times P \times Q,$$

wobei $P = \mathbb{P}(\mathbb{C}^{N*})$, $Q = \mathbb{P}(\mathbb{C}^N)$ und $T_x$ der Tangentialraum der Faser im Punkt $x$ sei. Die Dimension von $C(X/D)$ ist $N$. Man erhält ein kommutatives Diagramm:

---

[12] Vgl. Fußnote 8.
[13] Zusammen mit dem Beweis in [L2] folgt möglicherweise auch die starke simultane Auflösbarkeit.

$$B_Y C(X/D) \xrightarrow{\kappa'} B_Y X$$
$$\downarrow b' \qquad \downarrow b$$
$$C(X/D) \xrightarrow{\kappa} X$$

Mit $\tilde{\kappa} = \kappa b'$ bezeichnen wir die Komposition. Bei den folgenden Aussagen interessieren wir uns nur für den Keim von $X$ über $(D, 0)$ und obere Schranken für die Dimension der zentralen Faser.

**A.1.1:** $B_Y C(X/D)$ hat Faserdimension $\leq N - 2$ genau dann, wenn $C(X/D)$ Faserdimension $\leq N - 2$ hat und für den Schnitt $X_H = X \cap H$ mit allgemeinen Hyperebenen $H = H' \times D$ der Raum $B_Y C(X_H/D)$ Faserdimension $\leq N - 3$ hat. In diesem Fall haben die Fasern von $X$ konstante Multiplizität längs $Y$.

Beweis: Für allgemeine Hyperebenen enthält $X \times P \times \mathbb{P}(H')$ keine Komponente des Urbilds von $C$ in $B_Y C(X/D)$. Der Durchschnitt $\tilde{\kappa}^{-1}(C) \cap X \times P \times \mathbb{P}(H')$ ist also von Dimension $< N - 1$. Dann ist der Durchschnitt $B_Y C(X/D) \cap X \times P \times \mathbb{P}(H')$ der Abschluss von $B_Y C(X/D) \cap X \times P \times \mathbb{P}(H')$ über $(X - C) \cap H$, d.h. die strikte Transformierte. Denn dieser Durchschnitt hat die reine Dimension $N - 1$ und kann keine Komponente über $C$ haben.

Sei $A \subseteq P \times Q$ die spezielle Faser von $B_Y C(X/D)$. $A$ ist auch die spezielle Faser des exzeptionellen Divisors $\tilde{E} = \tilde{\kappa}^{-1}(Y)$, der rein von Dimension $N - 1$ ist. Die vertikalen Komponenten von $\tilde{E}$ haben Dimension $N - 1$, die übrigen Komponenten von $A$ haben Dimension $N - 2$.

Sei $\dim A < N - 1$. Man kann dann $H'$ außerhalb des Bildes der Projektion $A \to P$ wählen. Wenn außerdem $B_Y C(X/D) \cap X \times P \times \mathbb{P}(H')$ wie oben eine strikte Transformierte ist, dann ist (für einen genügend kleinen Repräsentanten) die Abbildung

(*)  $B_Y C(X/D) \cap X \times P \times \mathbb{P}(H') \to B_Y C(X_H/D)$, $(x, h, l) \mapsto (x, h \cap H', l)$,

wohldefiniert, endlich und surjektiv. Denn das Urbild von $(x, h', l)$ ist analytisch und kompakt in $\{h | h \cap H' = h', h \nsubseteq H'\} \cong \mathbb{C}$. Wählt man $H'$ noch so, dass keine Komponente von $A$ in $P \times \mathbb{P}(H')$ enthalten ist, so hat $B_Y C(X_H/D)$ Faserdimension $< \dim A \leq N - 2$. Ferner folgt aus dem Diagramm: $C(X/D)$ hat Faserdimension $< N - 1$ und die Aufblasung $B_Y X$ hat Faserdimension $< 2$. Denn die exzeptionelle Menge $E$ von $B_Y X$ kann keine vertikale Komponente haben, da sie das Bild von $\tilde{E}$ ist. Also ist faserweise die Multiplizität konstant längs $Y$. (Vgl. [T, Ch. I, 5.1.1].)

Sei $\dim A = N - 1$. Das Bild von $A \to P$ ist auch die spezielle Faser von $C(X/D)$, bei surjektiver Abbildung ist sie von Dimension $N - 1$. Im Fall, dass $A \to P$ nicht surjektiv ist, wählt man wie oben $H'$ außerhalb des Bildes und mit $\dim A \cap P \times \mathbb{P}(H') = N - 2$. Dann folgt aus der Endlichkeit von (*), dass $B_Y C(X_H/D)$ Faserdimension $N - 2$ hat.

Wir betrachten jetzt den Durchschnitt $X_H = X \cap H$ mit $H = H' \times D$. Für allgemeine $H$ ist $(\{0\} \times D, (0,0))$ die kritische Menge der Projektion $(X_H, (0,0)) \to (D, 0)$. Es genügt dafür, dass $\bigl(C \cap H, (0,0)\bigr) = (Y, (0,0))$ und $H'$ nicht in der speziellen Faser von $C(X/D)$ liegt.

**A.1.2:** Für die Kurvenfamilie $X_H \subseteq H' \times D$ gilt: $B_Y C(X_H/D)$ hat Faserdimension $\leq N - 3$ genau dann, wenn dies für $C(X_H/D)$ gilt und $B_Y X_H$ endliche Fasern hat (d.h. die Multiplizität konstant ist).

Beweis: ⇒ folgt wie in A.1.1. ⇐: Mit den Bezeichnungen $P' = \mathbb{P}(H'^*)$, $Q' = \mathbb{P}(H')$ sei $A \subseteq P' \times Q'$ die spezielle Faser von $B_Y C(X_H/D)$. Wenn die Dimension $N-2$ ist und $A \to P'$ surjektiv ist, hat auch $C(X_H/D)$ diese Faserdimension. Anderenfalls hat $A \to P'$ nicht endliche Fasern und das Bild von $A \to Q'$ ist mindestens eindimensional. Dies ist auch die Faser von $B_Y X_H$.

Man kann auch die relative Nash-Transformation $N(X/D)$ betrachten. Diese trägt ein Vektorraumbündel $T \subseteq \mathbb{C}^N_{N(X/D)}$, die Fortsetzung des Tangentialbündels. Man hat eine surjektive Abbildung $\mathbb{P}((\mathbb{C}^N_{N(X/D)}/T)^*) \to C(X/D)$ durch Vernachlässigung der Basispunkte. Die Faserdimension auf der linken Seite ist die von $N(X/D)$ vermehrt um $N-3$. Die Bedingung, dass $N(X/D)$ höchstens eindimensionale Fasern hat, ist i.a. stärker als dass $C(X/D)$ Fasern der Dimension $\leq N-2$ hat.

Für eindimensionale Familien unterscheiden sich die Bedingungen über die Faserdimension der Nash-Transformation und des Konormalenraums aber nicht.

**A.1.3:** Die relative Nash-Transformation $N(X_H/D)$ ist endlich über $X_H$ genau dann, wenn der Konormalenraum $C(X_H/D)$ Faserdimension $\leq N-3$ hat.

Beweis: Sei $A \subseteq \mathbb{P}(H'^*)$ die spezielle Faser von $C(X_H/D)$, $B \subseteq \mathbb{P}(H')$ die spezielle Faser von $N(X_H/D)$. Dann gilt

$$A = \bigcup_{l \in B} \mathbb{P}((H'/l)^*) \subseteq \mathbb{P}(H'^*).$$

Die algebraische Menge ist daher eine endliche Vereinigung von Hyperebenen oder der ganze projektive Raum. Wenn $B$ endlich ist, ist also $A$ von Dimension $\leq N-3$. Die Umkehrung gilt auch, weil die Hyperebene $\mathbb{P}((H'/l)^*)$ die Gerade $l$ eindeutig bestimmt.

Wir nehmen nun zur Vereinfachung an, dass die kritische Menge $C = Y$ ist. Man beachte, dass dies auch die kritische Menge von $X_H$ ist, wenn $H'$ nicht in der speziellen Faser von $C(X/D)$ liegt.

**A.1.4:** Die Bedingungen in A.1.1 beschreiben die Whitney-Regularität von $(X-Y, Y)$. Analog beschreiben die Bedingungen in A.1.2 die Whitney-Regularität von $(X_H - Y, Y)$, wenn die kritische Menge $C = Y$ ist.

Beweis: Sei $f = (f_1, \ldots, f_m)$ ein Erzeugendensystem für das Ideal von $X$. Nach dem idealtheoretischen Bertini-Satz ([T, Ch. II, th. 2.1.1] in der Version von [G, th. 2.7]) gilt an Stellen $x \in Y$, $x \neq 0$ nahe $0$

(*) $\quad\quad\quad\quad f_t \in \overline{(x_1, \ldots, x_N)(f_{x_1}, \ldots, f_{x_N})} \subseteq \mathcal{O}^m_{X,x}.$

An nichtsingulären Stellen $x \notin Y$ der Fasern gilt sogar $f_t \in (f_{x_1}, \ldots, f_{x_N}) \subseteq \mathcal{O}^m_{X,x}$. Die Aufblasung $B_Y C(X/D)$ ist das projektive Schema zu dem Rees-Ring des Untermoduls $(x_1, \ldots, x_N)(f_{x_1}, \ldots, f_{x_N}) \subseteq \mathcal{O}^m_{X,x}$. Wenn die Faserdimension hiervon $\leq N-2$ ist, ist das Spezialisierungsprinzip für ganze Abhängigkeit von Moduln anwendbar und liefert (*) am Nullpunkt. (Man kann [KT2, th. (1)(b)] auf den ganzen Abschluss von $(x_1, \ldots, x_N)(f_{x_1}, \ldots, f_{x_N})$ und dessen Erweiterung um $f_t$ anwenden, vgl. [GK]. Ein direktes Argument geht so: Der Kegel zu dem normalisierten Rees-Ring von $(x_1, \ldots, x_N)(f_{x_1}, \ldots, f_{x_N})$ ist reindimensional und die spezielle Faser hat Kodimension mindestens 2. Daher ist die durch $f_t$ als Element vom Grad 1 definierte rationale Funktion über die spezielle Faser fortsetzbar.) Die Bedingung (*) für $x = 0$ ist gleichbedeutend mit der Whitney-Regularität in der Nähe des Nullpunkts [G, th. 2.5].

Wenn die Whitney-Regularität erfüllt ist, ist $V := \mathbb{C}^N \times \{0\}$ nicht in der speziellen Faser des (absoluten) Konormalenraums $C(X)$ enthalten. Daher ist die Abbildung,

$$B_Y C(X) \to B_Y C(X/D), (x, h, l) \mapsto (x, h \cap V, l),$$

wohldefiniert und surjektiv. Es ist bekannt [T, Ch. V, th. 1.2], dass $B_Y C(X)$ Faserdimension $\leq N - 2$ hat, also auch $B_Y C(X/D)$. Ein analoger Beweis gilt für $X_H$.

**Anhang 2: Darstellung der Nash-Transformation**

Wir beschreiben hier die bekannte Darstellung der Nash-Transformation als Aufblasung eines Ideals (vgl. z.B. [T]). Sei zunächst $(X, 0) \subseteq (\mathbb{C}^N, 0)$ eine einzelne zweidimensionale normale Singularität mit Ideal $I \subseteq R := \mathcal{O}_{\mathbb{C}^N, 0}$. Dann ist $I$ ein Primideal und $R_I$ regulär. Also kann man Parameterelemente $g_1, \ldots, g_{N-2} \in I$ wählen, die $IR_I$ erzeugen. Dann ist $V(g_1, \ldots, g_{N-2})$ ein vollständiger Durchschnitt und es gibt ein $h \notin I$ mit $(g_1, \ldots, g_{N-2})R_h = IR_h$. Also gilt

$$V(g_1, \ldots, g_{N-2}) - V(h) = X - V(h), X \nsubseteq V(h).$$

Sei $J := I_{N-2}(g_{x_1}, \ldots, g_{x_N})\mathcal{O}_X$ das Ideal der $N - 2$−Minoren der Jacobi-Matrix $(g_{x_1}, \ldots, g_{x_N})$ von $g = (g_1, \ldots, g_{N-2})$. Dieses Ideal ist inversibel außerhalb der Singularität 0. Denn sei $x \in X - \{0\}$ und $\varphi_1, \ldots, \varphi_{N-2}$ ein Erzeugendensystem von $I$ nahe $x$. Wegen $(g_1, \ldots, g_{N-2}) \subseteq I$ gilt

$$g_i = a_{i1}\varphi_1 + \cdots + a_{i,N-2}\varphi_{N-2}, i = 1, \ldots, N - 2, a_{ij} \in \mathcal{O}_{\mathbb{C}^N, x}.$$

Für die Ableitungen folgt

$$g_{i, x_k} \equiv a_{i1}\varphi_{1, x_k} + \cdots + a_{i, N-2}\varphi_{N-2, x_k} \mod(\varphi_1, \ldots, \varphi_{N-2}), k = 1, \ldots, N.$$

Weil die Matrix $(\varphi_{j, x_k})_{\substack{1 \leq j \leq N-2 \\ 1 \leq k \leq N}}$ einen invertierbaren $N - 2$−Minor enthält, folgt daraus

$$I_{N-2}(g_{x_1}, \ldots, g_{x_N})\mathcal{O}_{X, x} = \det(a_{ij})\mathcal{O}_{X, x}.$$

Die Aufblasung von $J$ ist die Nash-Transformation. Denn nach Definition ist die Nash-Transformation $\nu: N(X) \to X$ der Abschluss des Bildes der Abbildung in das Grassmann-Bündel

$$X_{reg} \to X \times G, x \mapsto (x, T_x),$$

die regulären Punkten ihren Tangentialraum zuordnet, zusammen mit der ersten Projektion. Aus dem universellen Bündel auf $G$ erhält man ein Bündel $T \to N(X)$, welches das Tangentialbündel von $X_{reg}$ auf $N(X)$ fortsetzt. Mittels der Plücker-Einbettung kann $N(X)$ auch als Abschluss des Bildes von

$$X_{reg} \to X \times \mathbb{P}(\Lambda^{N-2}(\mathbb{C}^N)^*), x \mapsto (x, (\mathbb{C}^N/T_x)^*) \mapsto (x, \Lambda^{N-2}(\mathbb{C}^N/T_x)^*),$$

angesehen werden. Wenn $\delta_1, \ldots, \delta_m$ die $m = \binom{N}{N-2}$ maximalen Minoren der Jacobi-Matrix $(g_{x_1}, \ldots, g_{x_N})$ sind, so ist diese Abbildung für $x \in X_{reg} - V(h)$ darstellbar als

$$x \mapsto [g_{x_1}(x), \ldots, g_{x_N}(x)] \mapsto [\delta_1(x), \ldots, \delta_m(x)],$$

und der Abschluss ist gerade die Aufblasung von $J = (\delta_1, \ldots, \delta_m)$. Das Ideal $J\mathcal{O}_{N(X)}$ ist die Garbe der Schnitte des dualen Universalbündels auf $N(X) \subseteq X \times \mathbb{P}(\Lambda^{N-2}(\mathbb{C}^N)^*)$ und damit des Geradenbündels

$$L := \Lambda^{N-2}(\mathbb{C}^N/T)^{**} = \Lambda^{N-2}(\mathbb{C}^N/T) \cong \Lambda^2 T^*.$$

Wir bemerken, dass $N(X) \to X$ nur dann endlich ist, wenn $X$ glatt ist. Denn wegen der Normalität folgt dann $N(X) \cong X$ und $X$ hat genau einen Limes $T_0$ von Tangentialräumen bei 0. Sei $q: \mathbb{C}^N \to \mathbb{C}^2$ eine lineare Projektion, so dass $X \to \mathbb{C}^2$ endlich und $T_0 \to \mathbb{C}^2$ injektiv ist. Wenn $X$ nicht glatt ist, gibt es eine kritische Kurve in $X$, wo die Tangentialräume vertikal sind, d.h. nicht injektiv abgebildet werden. Dann gibt es aber auch einen vertikalen Limes.

Sei $(X, 0) \subseteq (\mathbb{C}^N \times \mathbb{C}, 0)$ nun eine Deformation von zweidimensionalen normalen Singularitäten, d.h. die kritische Menge $(C, 0)$ der Projektion $p: (X, 0) \to (\mathbb{C}, 0)$ auf die zweite Komponente ist endlich über $(\mathbb{C}, 0)$. Dann kann man einen vollständigen Durchschnitt $V(g_1, \ldots, g_{N-2}) \supseteq X$ wählen, der außerhalb einer Hyperfläche $V(h)$, die die Faser $X_0$ nicht enthält, mit $X$ übereinstimmt. Es genügt, entsprechende Elemente für die Faser innerhalb des Ideals von $X$ beliebig fortzusetzen. Insbesondere ist das Minorenideal der relativen Jacobi-Matrix $J := I_{N-2}(g_{x_1}, \ldots, g_{x_N})\mathcal{O}_X$ inversibel außerhalb von $C$ und die Einschränkung auf $X_0$ ist nicht Null.

Die Aufblasung von $J$ ist nun die relative Nash-Transformation $\nu: N(X/D) \to X$, also der Abschluss des Bildes der Abbildung in das Grassmann-Bündel

$$X - C \to X \times G, x \mapsto (x, T_x),$$

die glatten Punkten der Fasern den Tangentialraum der Faser zuordnet. Wieder gibt es eine Fortsetzung $T \to N(X/D)$ des relativen Tangentialbündels und $J\mathcal{O}_{N(X/D)}$ ist die Garbe der Schnitte von $\Lambda^2 T^*$.

**Literaturangaben**